\documentclass [10pt,reqno]{amsart}
\usepackage {hannah}
\usepackage {color,graphicx}
\usepackage{texdraw}
\usepackage{epsfig}
\usepackage{epstopdf}
\usepackage [latin1]{inputenc}
\usepackage[all]{xy}

\newcommand{\RR}{\mathbb R}

\newcommand{\ZZ}{\mathbb Z}

\newcommand{\TP}{\mathbb{TP}}
\newcommand{\cT}{\mathcal T}
\newcommand{\cG}{\mathcal G}

\DeclareMathOperator {\ev}{ev}

\DeclareMathOperator {\lab}{lab}
\DeclareMathOperator {\calM}{\mathcal{M}}

\newcommand {\dunion}{\,\mbox {\raisebox{0.25ex}{$\cdot$}  
\kern-1.83ex $\cup$} \,}

\title [The space of tropically collinear points is shellable]{The space of tropically collinear points is shellable}

\author {Hannah Markwig}
\address {University of Michigan,
Department of Mathematics,
2074 East Hall,
530 Church Street,
Ann Arbor,
MI 48109-1043}
\email{markwig@umich.edu}

\author{Josephine Yu}
\address{Massachusetts Institute of Technology,
    Department of Mathematics,
    Cambridge, MA 02139-4307}
\email {jyu@math.mit.edu}

\begin {document}

 \begin{abstract}
The space $T_{d,n}$ of $n$ tropically collinear points in a fixed tropical projective space $\mathbb{TP}^{d-1}$ is equivalent to the tropicalization of the determinantal variety of matrices of rank at most $2$, which consists of real $d\times n$ matrices of tropical or Kapranov rank at most $2$, modulo projective equivalence of columns.  We show that it is equal to the image of the moduli space $\mathcal{M}_{0,n}(\TP^{d-1},1)$ of $n$-marked tropical lines in $\TP^{d-1}$ under the evaluation map.  Thus we derive a natural simplicial fan structure for $T_{d,n}$ using a simplicial fan structure of $\mathcal{M}_{0,n}(\TP^{d-1},1)$ which coincides with that of the space of phylogenetic trees on $d+n$ taxa.   The space of phylogenetic trees has been shown to be shellable by Trappmann and Ziegler.  Using a similar method, we show that $T_{d,n}$ is shellable with our simplicial fan structure and compute the homology of the link of the origin.  The shellability of $T_{d,n}$ has been conjectured by Develin in \cite{Dev05}.
 \end {abstract}

\maketitle

\section{Introduction}

Let $(\RR, \oplus, \odot)$ be the {\em tropical semiring} where the tropical addition $\oplus$ is taking minimum and the tropical multiplication $\odot$ is the usual addition.  We will work in the tropical projective space $\TP^{d-1} = \RR^d / (1,\dots, 1)\RR$ obtained by quotienting out the tropical scalar multiplication (tropical projective equivalence).

A {\em tropical line} in $\TP^{d-1}$ is a one dimensional polyhedral complex in $\TP^{d-1}$ which is combinatorially a tree with unbounded edges in directions $e_1, \dots, e_d$ and the {\em balancing condition} at each vertex as follows.  At a vertex $V$, let $u_1, \dots, u_k$ be the primitive integer vectors pointing from $V$ to its adjacent vertices (respectively, in direction of the unbounded edges adjacent to $V$).  The balancing condition holds at $V$ if  $u_1 + \cdots + u_k = 0$ in $\TP^{d-1}$.
A configuration of $n$ points in $\TP^{d-1}$ is called \emph{tropically collinear} if there is a tropical line which passes through the $n$ points.  Let $T_{d,n}$ be the space of all such configurations.  An element of $T_{d,n}$ is represented by a real $d \times n$ matrix whose $n$ columns are representatives of the $n$ points in $\TP^{d-1}$.  

The columns of a matrix are tropically collinear if and only if the matrix lies in the tropical variety of the determinantal ideal generated by  $3 \times 3$ minors of a $d\times n$ matrix of indeterminates, which is a polyhedral subfan of the Gr\"{o}bner fan consisting of those cones whose corresponding initial ideal is monomial-free.  Since $T_{d,n}$ is equal to this tropical variety modulo projective equivalence of the columns in the matrix, it is also a polyhedral fan.  However this fan structure is not simplicial, and the topology of the link of the origin is difficult to analyze.  Develin conjectured in \cite{Dev05} that $T_{d,n}$ is shellable for all $d$ and $n$, and proved his conjecture for $d=3$ (or $n=3$). In this paper we prove his conjecture for all $d$ and $n$, with a refined fan structure.  We give a triangulation of the fan that lets us treat $T_{d,n}$ as a subcomplex of the space of phylogenetic trees, which in turns gives us a way to prove strong combinatorial properties such as shellability.  Shellability implies that our space is Cohen-Macaulay and has homology only in the top dimension.  Moreover, a shelling order gives us a way to compute the top homology of the link of the origin.  

In Section 2, we will derive a simplicial fan structure on $T_{d,n}$, using moduli spaces of tropical curves and the space of phylogenetic trees $\cT_{n+d}$.  
Since $T_{d,n}$ is closed under simultaneous translation of all points
, we will mod out by this action and obtain a pointed simplicial fan that we denote by $T_{d,n}'$.   We then intersect this fan with the unit sphere centered at the origin to obtain a simplicial complex which we will call $T_{d,n}''$.  The simplicial complex structures of $T_{d,n}$, $T_{d,n}'$, and $T_{d,n}''$ are all the same.  In particular, if one of them is shellable, then so are the other two.

 A parametrized tropical line can be thought of as an abstract tropical curve (a leaf-labeled tree) $\Gamma$ together with a map \begin{displaymath}h:\Gamma\rightarrow \TP^{d-1},\end{displaymath} such that the image $h(\Gamma)$ is a tropical line as defined above. Our parametrized tropical lines are equipped with certain marked points $x_i$. In section $2$, we will recall the definition of moduli spaces of $n$-marked parametrized tropical lines, and evaluation maps which send a tuple $(\Gamma,h,x_i)$ to $h(x_i)\in \TP^{d-1}$. 
 We will show that $T_{d,n}$ is the image of the moduli space of $n$-marked parametrized tropical lines under the evaluation map. 

Moduli spaces of tropical curves can be used to derive results in enumerative tropical geometry. This is why these moduli spaces attracted a lot of attention recently (see e.g.\ \cite{Mi07}, \cite{GM053} or \cite{GKM07}). Their simplicial fan structure equals the structure of the space of trees, $\cT_{n+d}$ (see \cite{GKM07}). In fact, we can identify $T_{d,n}$ with a subcomplex of the space of trees $\cT_{n+d}$ on which the evaluation map is injective, the subcomplex induced on the vertices corresponding to ``bicolored splits'' (see definition \ref{def-bicoloredsplit}).

In \cite{TZ98}, Trappmann and Ziegler showed that the space of trees $\cT_{n+d}$ is shellable. 
 Since we derive our simplicial complex structure for $T_{d,n}'$ using the space of trees $\cT_{n+d}$, we will use a similar method to show that the space $T_{d,n}'$ is shellable, in Section 3.  We compute the homology of $T_{d,n}''$ in Section 4.
Our main results can be summarized as follows:

\begin{theorem}
The simplicial complex $T_{d,n}''$ is shellable and has the homotopy type of a wedge of $n+d-4$-dimensional spheres.  The number of spheres is equal to the number of 
simultaneous partitions of an $(n-1)$-set and a $(d-1)$-set into the same number of non-empty ordered parts.  This number equals 
$$\sum_{k=1}^{\min(n-1,d-1)} (k!)^2 S(n-1,k) S(d-1,k),$$ 
$$\text{ where } S(m,k) = \frac{1}{k!} \sum_{i=1}^k (-1)^{k-i} {k \choose i} i^m \text{ is the Stirling
number of the second kind.}$$
\end{theorem}
This theorem follows from the theorems \ref{thm-shell} and \ref{thm-hom}.


In general, not much is known about the algebraic topology of tropical varieties.  Not all tropical varieties have homology only in the top dimension \cite[Example 5.2]{ST07}.  However, there may be many classes of tropical varieties that have only top homology and are even shellable.
For example, it is not known whether the tropical varieties of generic complete intersections studied in \cite{ST07} are shellable.
Hacking \cite{Hac} gave some sufficient conditions for some tropical varieties to have only top homology.  He also gave a moduli space interpretation of $T_{d,n}$ and showed that it has only top homology.

{\bf Acknowledgments}.  We would like to thank Bernd Sturmfels for helpful conversations.  This work began when the authors were at the Institute for Mathematics and Its Application in Minneapolis.  Josephine Yu was supported by a Clay Liftoff Fellowship during summer 2007.

\section{The space $T_{d,n}'$ inside the space of phylogenetic trees}

As mentioned in the introduction, we want to derive a simplicial fan structure for $T_{d,n}$ using moduli spaces of tropical curves.
Let us start by defining tropical curves and their moduli spaces. 

\subsection{The space $\calM_{0,N}$ of $N$-marked abstract tropical curves}

An \emph{abstract tropical curve} is a tree
$\Gamma$ whose vertices have valence at least $3$.
The internal edges are equipped with a finite positive length.  The leaf edges are considered unbounded.
An $N$-marked abstract tropical curve is a tuple $(\Gamma,x_1,\ldots,x_N)$ where
$\Gamma$ is an abstract tropical curve and $x_1,\ldots,x_N$ are
distinct unbounded edges. (For more details, see \cite{GM053}, definition 2.2.)
The space $ \calM_{0,N} $ is defined to be the \emph{space of
all $N$-marked abstract tropical curves} with exactly $N$ leaves.
The following picture shows an example of a $5$-marked abstract tropical curve:
\begin{center}
\includegraphics{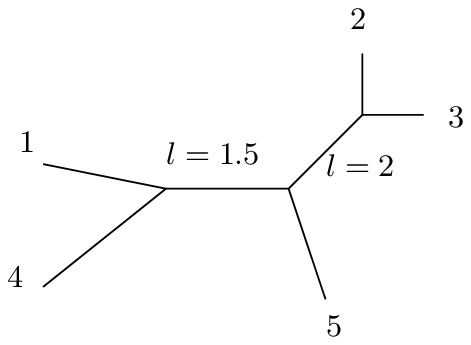}
\end{center}

Let $\cT_{N}$ be the space of phylogenetic trees on $N$ taxa. 
 A {\em phylogenetic tree} (or a semi-labeled tree or a leaf-labeled tree) on $N$ taxa is a tree with $N$ leaves labeled by $[N]$ and vertices of valence at least 3
 such that the internal edges have positive lengths and the leaf edges have non-negative lengths.  
In other words, a phylogenetic tree  on $N$ taxa is an $N$-marked abstract tropical  curve whose leaf-edges are assigned a non-negative length.
Hence $\calM_{0,N}$ is the space $\cT_N$ of phylogenetic trees modulo the space of {\em star trees}.  A star tree  is a phylogenetic tree with no internal edges. The following figure shows a star tree on  5 taxa.

\begin{center}
\includegraphics{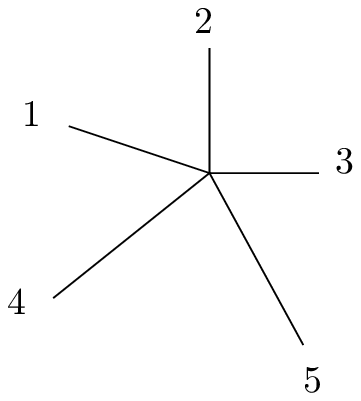}
\end{center}

\begin{proposition}\label{prop-m0nfan}
The spaces $\cT_N$  and $ \calM_{0,N} $ can be embedded as simplicial fans into real vector spaces of dimensions ${N \choose 2}$ and $\binom{N}{2}-N$ respectively. 
\end{proposition}

For a complete proof, see \cite[theorem 4.2]{SS04a} or \cite[theorem 3.4]{GKM07}. In fact, $ \calM_{0,N} $ is a tropical fan or a balanced fan, see \cite{GKM07, SS04a, ST07}.  These two fans are the fans $\cG_{2,N} \cap \RR_{+}^{N \choose 2}$ and $\cG''_{2,N}$ respectively in the tropical Grassmannian \cite{SS04a}.

As an idea why proposition \ref{prop-m0nfan} is true, note that a phylogenetic tree is completely determined by the metric $d$ it induces on the set $[N]$:  the distance $d_{ij}$ between two elements $i,j \in [N]$ is the sum of the edge lengths along the unique path between $i$ and $j$ in the tree.  Hence $\cT_N$ can be embedded in $\RR^{N \choose 2}$, and $\calM_{0,N}$ can be embedded in $\RR^{N \choose 2}$ modulo the $N$-dimensional vector space of star trees.  The simplicial complex of $\cT_N$ is a cone over the simplicial complex of $\calM_{0,N}$.
 
Let us now recall the embedding and the fan structure of $\cT_N$.
The rays of $\cT_{N}$ correspond to {\em splits}, partitions of $[N]$ into two non-empty parts.  Removing an edge in a phylogenetic tree decomposes the tree into two connected components, inducing a split on the set of leaf labels.  Two splits $A|A'$ and $B|B'$ are {\em compatible} if at least one of $A \cap B$, $A \cap B'$, $A' \cap B$, and $A' \cap B'$ is empty.  A set of splits is pairwise compatible if and only if there is a (unique) tree whose edges induce exactly those splits \cite[theorem 2.35]{PS05}.
We sometimes do not distinguish between a tree and its corresponding set of splits.  For a split $A|B$, we can write it as $A~|~([N]\setminus A)$  where  $[N]:=\{1,\ldots,N\}$, so we can respresent this  split by just $A$.
Edges of a tree are labelled with their split or with the part $A$ of the split not containing $1$, which is also called the {\em edge label}.
A set $S$ of splits generate a cone in $\cT_N$ if and only if the splits in $S$ are pairwise compatible.  The singleton splits $\{i\} | ( [N]\backslash\{i\}), i \in N$ are compatible with all other splits, so every maximal cone in $\cT_N$ contains the cone generated by those splits.


\subsection{The space $\calM_{0,n}(\TP^{d-1},1)$ of $n$-marked tropical lines in $\TP^{d-1}$}

Now we will review how the $N$-marked abstract tropical curves parametrize tropical lines.

\begin{definition}
Let $N=n+d$. A \emph {(parametrized) $n$-marked tropical line in} $ \TP^{d-1} $ is a
  tuple $ (\Gamma,x_1,\dots,x_N,h) $, where $
  (\Gamma,x_1,\dots,x_N) $ is an abstract $N$-marked tropical curve and $ h:
  \Gamma \to \TP^{d-1}  $ is a continuous map satisfying:
  \begin {enumerate}
  \item On each edge of $ \Gamma $ the map $h$ is of the form $ h(t) = a + t
    \cdot v $ for some $ a \in \RR^r $ and $ v \in \ZZ^r $. The integral
    vector $v$ occurring in this equation if we parametrize $E$ by an interval $[0,l(E)]$ (starting at $V\in \partial E$) will be denoted $ v(E,V) $ and called the \emph {direction} of $E$ at $V$. If $E$ is an unbounded edge and $V$ its only boundary point we will write for simplicity $ v(E) $ instead of $ v(E,V) $.
  \item For every vertex $V$ of $ \Gamma $ we have the \emph {balancing
    condition}
      \[ \sum_{E| V \in \partial E} v(E,V) = 0. \]
  \item $ v(x_i)=0 $ for $ i=1,\dots,n $ --- i.e.\ each of the first $n$ ends is
    contracted by $h$.
    \item $ v(x_i) =e_{i-n}$ for $ i>n $ --- i.e.\ the
    remaining $ N-n $ ends are mapped to the $d$ canonical directions of $ \TP^{d-1} $.
  \end {enumerate}
  
  We will call the contracted ends the ``marked ends" and the directed unbounded ends the ``unmarked ends".  Note that this notion of markedness is different from the marked ends of an abstract tropical curve.
The space of all labeled $n$-marked tropical lines in $\TP^{d-1} $ will be denoted $ \calM_{0,n}(\TP^{d-1},1) $. 
\end{definition}

\begin{remark}\label{rem-labelends}
Note that this definition is a special case of definition 4.1 of \cite{GKM07}. As we are working with lines, the unmarked ends are mapped to different directions. Hence they are distinguishable by their direction and we do not need to label them as in \cite{GKM07}. In our case, there is no such difference between $\calM^{\lab}_{0,n}(\RR^{d},\Delta) $ and $ \calM_{0,n}(\RR^d,\Delta) $ as mentioned in construction 4.3 of \cite{GKM07}.
To keep notations simple, we will still label the ``unmarked ends'' according to the following rule: the end with direction $e_i$ gets the label $n+i$.
\qed
\end{remark}

The following picture shows an example of an element of $ \calM_{0,2}(\TP^{2},1) $.
We will always draw the marked (contracted) ends with a dotted line and the unmarked (directed) ends solid.
\begin{center}
\includegraphics{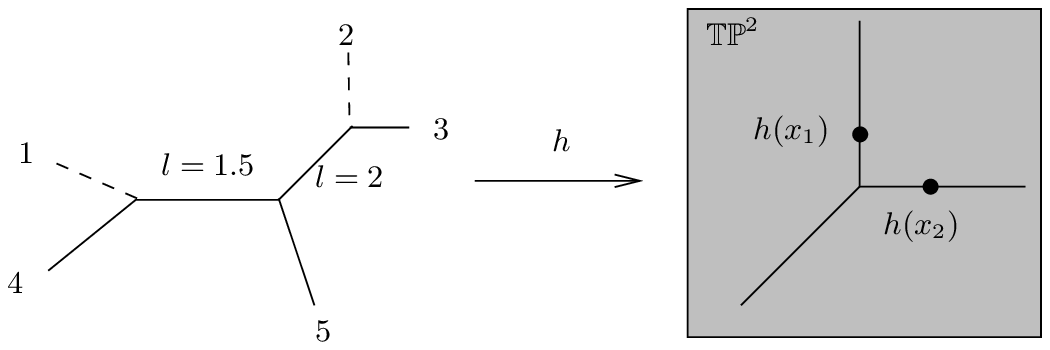}
\end{center}

Notice also that the image of a parametrized $0$-marked tropical line is a tropical line as defined in the introduction. We only need to check that all direction vectors are primitive integral vectors. But this is true because all direction vectors are sums of different canonical vectors, $v=e_{i_1}+\cdots+e_{i_r}$.  This observation also shows that the internal edges cannot be contracted, i.e.\ the direction vector $v$ of an internal edge cannot be zero.  Hence the space of tropical lines in $\TP^{d-1}$ is $\calM_{0,0}(\TP^{d-1},1)$. 

\begin{definition}
  For each $ i=1,\dots,n $ define the \emph{$i$-th evaluation map} $\ev_i$ by
  \begin{eqnarray*}
    \ev_i:
      \calM_{0,n}(\TP^{d-1},1)&\rightarrow&\TP^{d-1}\\
      (\Gamma,x_1,\dots x_N,h) & \longmapsto & h(x_i)
  \end{eqnarray*}
This is well-defined
  for the contracted ends since for them $ h(x_i) $ is a point in $ \TP^{d-1} $. The product $\ev=\ev_1\times \ldots\times \ev_n:\calM_{0,n}(\TP^{d-1},1)\rightarrow (\TP^{d-1})^n$ is called the \emph{evaluation map}.
  \end{definition}
 One can see easily that evaluation maps are linear on each cone of the fan $\calM_{0,n}(\TP^{d-1},1)$ (see e.g.\ example 3.3 of \cite{GM053}). In fact, it is even a tropical morphism (for more details, see \cite{GKM07}).
\begin{lemma}\label{lem-collinear}
The following are equivalent for a real $d\times n$ matrix $M$:
\begin{enumerate}
\item the columns of $M$ are tropically collinear in $\TP^{d-1}$,
\item the rows of $M$ are tropically collinear in $\TP^{n-1}$,
\item $M$ has Kapranov rank at most $2$,
\item $M$ has tropical rank at most $2$,
\item the columns of $M$ are in the image of the map $$\ev : \calM_{0,n}(\TP^{d-1},1) \rightarrow (\TP^{d-1})^n,$$
\item the rows of $M$ are in the image of the map $$\ev : \calM_{0,d}(\TP^{n-1},1) \rightarrow (\TP^{n-1})^d.$$
\end{enumerate}
\end{lemma}

\begin{proof}
The equivalence of (a),(b), and (c) is the definition of the Kapranov rank.  The equivalence of (c) and (d) follows from \cite[theorem 6.5]{DSS05}.
If the columns of $M$ are in the image of $\ev$ that means that they are $n$ distinguished points $(h(x_1),\ldots,h(x_n))\in (\TP^{d-1})^n$ on a tropical line $h(\Gamma)$. So obviously they are tropically collinear.
Given $n$ tropically collinear points $(p_1,\ldots,p_n)$, there is a tropical line $L$ through the points. 
We can find an abstract tropical curve $\Gamma$ and a map $h:\Gamma \rightarrow \TP^{d-1}$ parametrizing $L$. Then we attach new marked ends at the preimages of the $p_i$ and required those to be contracted by $h$. In this way, we construct a preimage of $(p_1,\ldots,p_n)$ under $\ev$.
This proves the equivalence of (a) and (e). The equivalence of (b) and (f) can be shown analogously.
\end{proof}
  
We will now give a simplicial fan structure for $\calM_{0,n}(\TP^{d-1},1)$. Define the \emph{forgetful map} $\Psi$ which forgets the map $h$ as
\begin{eqnarray*}
\Psi: \calM_{0,n}(\TP^{d-1},1)&\rightarrow &\calM_{0,n+d}\\
(\Gamma,x_1,\ldots,x_N,h)&\mapsto&(\Gamma, x_1,\ldots,x_N).
\end{eqnarray*}

The following proposition that we cite from \cite{GKM07} determines a fan structure of $\calM_{0,n}(\TP^{d-1},1)$, using the fact that $\calM_{0,n+d}$ is a simplicial fan by proposition \ref{prop-m0nfan}.
\begin{proposition}\label{prop-stablemapsfan}
The map 
\begin{eqnarray*}
\ev_1\times \Psi: \calM_{0,n}(\TP^{d-1},1)&\rightarrow & \TP^{d-1}\times\calM_{0,n+d}\\
(\Gamma,x_1,\ldots,x_N,h)&\mapsto& (\ev_1(\Gamma,x_1,\ldots,x_N,h),\Psi(\Gamma, x_1,\ldots,x_N,h))
\end{eqnarray*}
is a bijection.
\end{proposition}
For a proof, see \cite[proposition 4.7]{GKM07}.
The idea why this is true is that we can deduce the direction vectors of all edges from the direction vectors which are prescribed for the unmarked ends (see lemma 4.6 of \cite{GKM07}).
Once the image $h(V)$ of one vertex is given --- in our case we choose the vertex of the marked end $x_1$ --- the map $h$ is completely determined by the direction vectors of all edges and their lengths, hence by the underlying abstract tropical curve.

\subsection{The space $T_{d,n}'$ as a subcomplex of $\calM_{0,n+d}$}
As a consequence of lemma \ref{lem-collinear}, we want to describe the space $T_{d,n}$ of $n$ tropically collinear points as the image of $\calM_{0,n}(\TP^{d-1},1)$ under $\ev$. 

As before, we let $N = n+d$, 
and think of a tree $T$ with $N$ leaves as an abstract tropical curve with $n$ marked ends and $d$ unmarked ends, where the unmarked end with the label $n+i$ gets the direction $e_i$ as in remark \ref{rem-labelends} above. 
For the ends, we will sometimes call the property of being marked or unmarked (contracted or directed) the ``color'' of the leaf.  
\begin{definition}
Define a map $\pi:\cT_N\rightarrow (\TP^{d-1})^n$ as follows.
First we define it for splits. Let $S = A|B$ be a non-singleton split with $1 \in A$.  Let $u \in \RR^d$ be the sum of $e_i$ such that $i+n \in B$.  Let $\pi(S)$ be the $d\times n$ matrix whose $j^{\text{th}}$-column is $0$ if $j \in A$ and $u$ otherwise.  For a singleton split $S$, we define $\pi(S)=0$.
We extend $\pi$ linearly on each cone of $\cT_{d,n}$.
\end{definition}
Note that $\pi$ is $0$ on the star tree because star tree contain only the singleton splits.

\begin{lemma}\label{lem-commute}
The following diagram is commutative:

 \[ \begin{xy}
       \xymatrix{ 
\calM_{0,n}(\TP^{d-1},1) \ar[d]^{\Psi}\ar[r]^{\;\;\;\;\;\;\;\;\;\;\ev} & T_{d,n}\ar[d]\\
 \calM_{0,n+d} \ar[r]^{\pi}& T_{d,n}' & }
\end{xy} \]
where the vertical arrow on the right is modding out by translation of the point configuration in $\TP^{d-1}$.
\end{lemma}

\begin{proof}
We only have to check that $\pi$ maps a tree corresponding to a split $S=A|B$ to a tuple in $(\TP^{d-1})^{n-1}$ consisting of the positions of the images of the marked points relative to $h(x_1)$.
Since both $\pi$ and $\ev$ are linear on a cone, the commutativity of the diagram above follows.
To see this, note that $ S=A|B$ corresponds to a tree with exactly one bounded edge of length 1 such that the ends marked by the numbers in $A$ are on one side and the ends marked by $B$ on the other.
The following picture shows an example where $n=5$, $d=3$, $A=\{1,3,6,7\}$ and $B=\{2,4,5,8\}$:
\begin{center}
\includegraphics{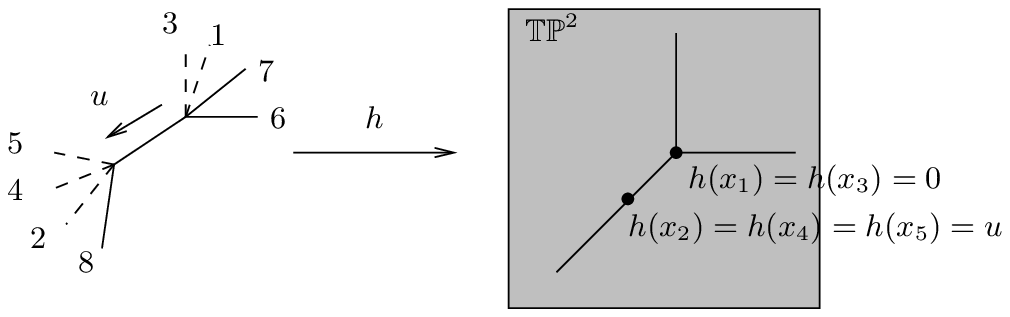}
\end{center}
Let us check the positions $h(x_i)$ of the marked ends relative to $h(x_1)$.
For all $i\in A, i\leq n$ (i.e.\ for all marked ends in $A$), $h(x_i)=h(x_1)$ so the relative evaluation is $0$. Let $u$ denote the direction vector of the bounded edge. By the balancing condition it is equal to the sum of all $e_i$ such that $i+n\in B$. Each marked point $x_j$ with $j\in B$ is mapped to $h(x_j)=h(x_1)+u$. The relative position to $h(x_1)$ is thus $u$.

This is precisely the definition of $\pi$.
\end{proof}

As a consequence of this lemma, we can think of $\pi$ as a ``relative'' or ``reduced'' version of the evaluation map $\ev$. 
Because of lemma \ref{lem-collinear}, we want to understand $T_{d,n}$ as the image of $ \calM_{0,n}(\TP^{d-1},1)$ under $\ev$. But since we mod out by simultaneous translations of all points in the definition of $T_{d,n}'$, we can do without the information $\ev_1(\Gamma,x_1,\ldots,x_N,h)=h(x_1)$ and consider the image of $\cT_{n+d}$ under $\pi$ instead.
We will make this more precise in lemma \ref{lem-pitree}.

\begin{definition}\label{def-bicoloredsplit}
Let $\cT_{d,n}$ be the subfan of $\cT_{N}$ consisting of trees whose non-singleton splits contain both marked and unmarked elements on each side.  We call such splits \emph{bicolored}. 
\end{definition}
 The subfan $\cT_{d,n}$ is an induced subcomplex of $\cT_N$ on the rays corresponding to those bicolored splits. 

Let $S=A|B$ be a split which is not bicolored. Let $1\in A$ and assume first that $A$ or $B$ contain only marked ends. Then the vector $u=\sum_{i:i+n\in B} e_i=0$ and $\pi(S)=0$. If $B$ contains no marked ends, $u$ is not necessarily $0$, but no column is equal to $u$. Hence $\pi(S)=0$ in this case, too. 
It follows immediatly that $\pi$ is not injective on a cone of $\cT_N$ if one of the generating rays corresponds to a non-bicolored split.
In fact, the subfan $\cT_{d,n}$ generated by rays corresponding to bicolored splits is precisely the union of all closed maximal cones of $\calM_{0,n}$ on which the relative evaluation map $\pi$ is injective. This is shown in the following lemma. 

\begin{lemma}\label{lem-pitree}
For any configuration of collinear points $p_1, \dots, p_n \in \TP^{d-1}$, there is a unique {\em canonical} tropical line through them with the property that if we attach a marked end at each point then we get a tree $T \in \cT_N$ with only bicolored splits.
Moreover, $\pi(T) = (p_1,\dots,p_n)$.
\end{lemma}

\begin{proof}
The canonical tropical line can be constructed as in \cite[Section 3]{Dev05}.  First take the tropical convex hull of the marked points $p_1, \dots, p_n$, which is the union of tropical line segments between pairs of points.  This is a bounded one-dimensional polyhedral complex which is combinatorially a tree.  Then there is a unique way to attach unbounded rays such that the balancing condition holds  \cite[Section 3]{Dev05}.  

An edge in a tropical line is called bounded if each end point is either one of the marked points or a vertex of the polyhedral complex.  We get a tree $T$ with bicolored splits after attaching marked ends at the marked points if and only if all bounded edges of the tropical line lie on a path between two marked points.  In other words, the union of bounded edges must be the tropical convex hull of $p_1, \dots, p_n$.  The canonical line is unique with this property.

The phylogenetic tree $T$ can be considered as an element of $\calM_{0,N}$, and the canonical line is an embedding of $T$.  Hence $T$ with this embedding is an element of $\calM_{0,n}(\TP^{d-1},1)$. By lemma \ref{lem-commute}, $\pi(T) = (p_1, \dots, p_n)$.

\end{proof}

\begin{corollary}
The map $\pi : \cT_{d,n} \rightarrow T_{d,n}'$ is an isomorphism of polyhedral fans, i.e.\ it is linear on each cone and is bijective.
\end{corollary}

We sum up the results of this section in the following proposition. We will use it to define a shelling order of $T_{d,n}'$ in the next section.
Recall that a subcomplex is {\em induced} means that a cell is in the subcomplex if and only if all its vertices are.

\begin{proposition}\label{prop-inducedsubcplx}
The space $T_{d,n}'$ has a triangulation that is isomorphic to the induced subcomplex of the space of phylogenetic trees $\cT_{d+n}$ on the vertices corresponding to bicolored splits, i.e.\ the splits containing both marked and unmarked leaves on each side.
\end{proposition}

\section{Shelling of $T_{d,n}'$}

In this section we will prove that $T_{d,n}'$ is shellable, in a similar way as in \cite{TZ98} where it was shown that the space of trees $\cT_N$ is shellable. We will use the description of $T_{d,n}'$ as the induced subcomplex of $\cT_N$ on the vertices of bicolored splits (see proposition \ref{prop-inducedsubcplx} and definition \ref{def-bicoloredsplit}).

In the following, we denote by $x\in C$ a vertex of a facet $C$ of a simplicial complex.
\begin{definition}
A {\em shelling} is an ordering of the facets of a pure-dimensional simplicial complex  
such that: For any two facets $C'<C$ there exist $C''$ and $x \in C$ such that
\begin{enumerate}
\item $C''<C$,
\item $x\notin C'$ and
\item $C\setminus x\subset C''$.
\end{enumerate}
\end{definition}
In \cite{TZ98} it is shown that these three conditions are equivalent to $C'\cap C\subset C''\cap C$, $C''<C$ and $C''$ differs from $C$ in only one element, $C\setminus C''=\{x\}$. The latter conditions are more common to define shellings.

Let us now recall the shelling of $\cT_N$ given in \cite{TZ98}.  First define an order on subsets of $[N]$ by:
\begin{displaymath}
A<B :\Leftrightarrow \max((A\setminus B)\cup (B\setminus A))\in B.
\end{displaymath}
To define the shelling order on the trivalent trees, we first ``split the trees along $1$'':
\begin{center}
\includegraphics{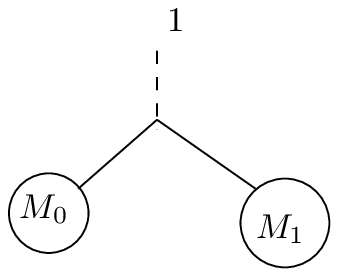}
\end{center}
There is a unique partition of $[N] \setminus \{1\}$ into two parts $M_0, M_1$ such that the two subtrees induced on the leaf-labels $M_0$ and $M_1$ are disjoint.  Let $M_0 < M_1$, i.e.\ $N \in M_1$.  Let $T_i$ be the subtree induced on $\{1\} \cup M_i$.

Let $T, T'$ be two trivalent trees, $\{T_0, T_1\}, \{T_0', T_1'\}$ be the pairs of subtrees, and $\{M_0, M_1\}$ and $\{M_0', M_1'\}$ be the corresponding pairs of leaf labels as above.  Then the order on the trees is defined recursively by:
$$
T' < T  ~\iff~ \left\{
\begin{array}{llll}
M_1' < M_1 & \text{or} &  &\\
M_1' = M_1 & \text{and} & T_1' < T_1 & \text{or}\\
T_1' = T_1 & \text{and} & T_0' < T_0
\end{array}
\right. 
$$

We can draw a tree respecting this order as in the figure above: then $1$ is the root of the tree, and the end vertex of $1$ is adjacent to two subtrees (with labels $M_0$ resp.\ $M_1$) that we order such that the bigger tree is on the right and the smaller on the left. Inductively, we can continue this process. According to this drawing of a tree, we will use the words ``child of a vertex $V$'' for a subtree that is below $V$.

We now come to the main result of the paper.
\begin{theorem}[Shelling of $T_{d,n}$]\label{thm-shell}
The complex $T_{d,n}$ is shellable.  The shelling order on $\cT_{N}$ defined above restricts to a shelling order on $T_{d,n}'$.
\end{theorem}

\begin{proof}
Suppose $T'<T$. We have to find $T''\in T_{d,n}$, $T'' < T$ such that $T''$ and $T$ differ only in one split $x \notin T'$. We use the identification of $T_{d,n}$ with the subcomplex induced on bicolored splits as described in proposition \ref{prop-inducedsubcplx}.

Without loss of generality we may assume that $M_1'<M_1$, for if $M_1=M_1'$, then there must be vertices $V$ of $T$ and $V'$ of $T'$ with the same leaf-labels below them such that the trees $T$ and $T'$ are equal above and to the right of $V$ and $V'$, and the right child of $V$ in $T$ is bigger than the right child of $V'$ in $T'$. Then the following argument works after replacing the end vertex of the marked end $1$ with $V$ or $V'$ and $N$ with the largest label below $V$ or $V'$.

Note that $M_1$ must have at least two leaves, since otherwise, the unique leaf would be labeled $N$, and there is no smaller $M_1'$ that also contains $N$.  Hence we can split $M_1$ into two subsets $L_0$ and $L_1$ satisfying $L_0<L_1$,  i.e.\ $N\in L_1$.

Let us first consider the case when $M_0$ contains at least two elements.  Since the split $M_0$ is bicolored, $M_0$ contains both marked and unmarked leaves.  In this case we swap subtrees of $T$ in the following way to obtain $T''$.

\begin{center}
\includegraphics{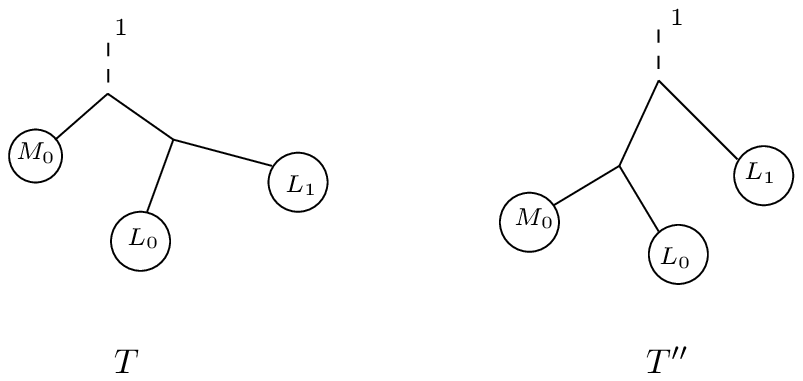}
\end{center}
We replace the split $M_1|\{1\}\cup M_0=L_0\cup L_1|\{1\} 
\cup M_0$ in $T$ by the split $L_0\cup M_0|\{1\}\cup L_1$ to obtain a new tree $T''$.  The split $\{1\}\cup M_0|L_0\cup L_1$ is not in the tree $T'$ because $M_1'<M_1$.
By our assumption  
we know that $L_0\cup M_0$ contains both marked and unmarked ends.  
As $1$ is marked and $N$ is unmarked and in $L_1$, $ \{1\}\cup L_1$ contains  both marked and unmarked ends, too.
As all other splits of $T''$ are splits of $T$, too, we conclude  
$T''\in T_{d,n}$.  Furthermore, $T''<T$, because $M_1''=L_1<L_0\cup L_1=M_1$.

Now let us consider the case when $M_0$ contains exactly one element.  Then $T$ has the form: 
\begin{center}
\includegraphics{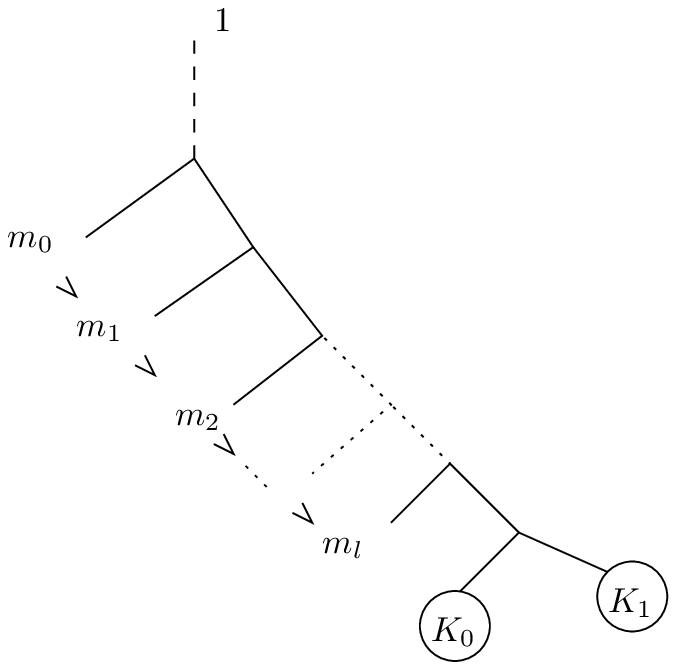}
\end{center}
where $m_0 > \cdots > m_l$ have the same color, and $K_0$ either contains a leaf of different color or a single leaf $m_{l+1} < m_l$ of the same color.   In the first case, we swap the split $K_0 \cup K_1$ with $K_0 \cup \{m_l\}$.  In the latter case, we swap $K_0 \cup K_1$ with $K_1 \cup \{m_l\}$.  This new tree $T''$ is smaller than $T$. Note that $l$ may be $0$.  However, we cannot have the following configuration because there is no smaller tree $T' < T$ that agrees with $T$ above and to the right of the vertex $V$:
\begin{center}
\includegraphics{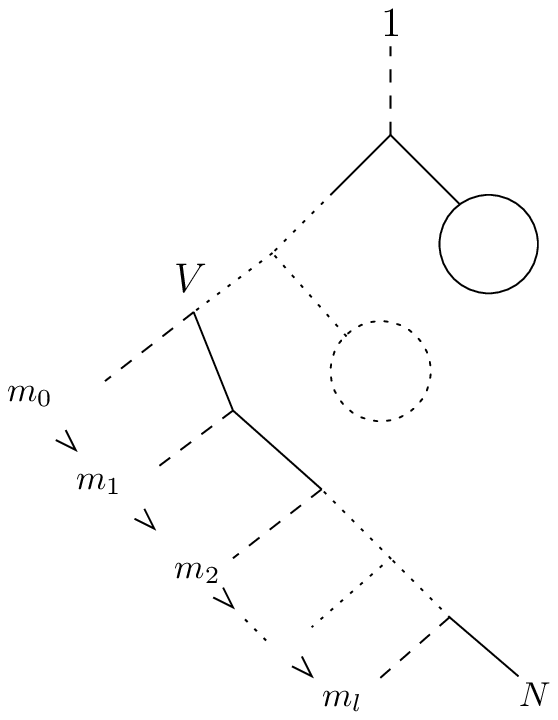}
\end{center}

We want to prove that $T'$ cannot contain the edge label $K_0 \cup K_1 
$. Assume it does. As $N\in K_1$, $K_0\cup K_1\subset M_1'$. Thus  
$M_0'$ can only consist of a subset of $\{m_0,\ldots,m_l\}$. But as  
all $m_i$ have the same color, $M_0'$ can contain at most one to be a bicolored split. Thus $M_1'=[N]\setminus\{1,m_i\}$ for some  
$i\neq 0$. But as $m_i<m_0$, $M_1'=[N]\setminus\{1,m_i\}>[N]\setminus  
\{1,m_0\}$ which is a contradiction. Hence $K_0\cup K_1\notin T'$ and  
we choose $x:=K_0\cup K_1$. 

\end{proof}

In \cite{Dev05}, it was shown that $T_{3,n}$ is shellable for all $n$ and a shelling order called ``snake ordering'' is given.  Our shelling order here is different from that.

\section{Homology of $T_{d,n}$}

In this section, we will use the shelling to compute the homology.
We have to count those trees $T\in T_{d,n}'$ which ``close a loop,''  
i.e.\ for each $x\in T$ there exists $T'\in T_{d,n}'$, $T'<T$ such  
that $T\setminus \{x\}\subset T'$.
In the proof of corollary 5 of \cite{TZ98}, it was shown that a tree with an internal left edge does not satisfy this condition, so all left edges must be leaves.  These types of trees are called {\em combs}.
\begin{center}
\includegraphics{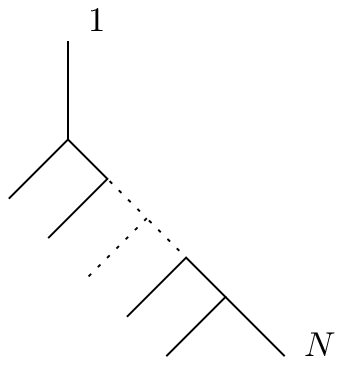}
\end{center}

Any comb such that the neighbour of $1$ is an unmarked end and the  
neighbour of $N$ is a marked end is in $T_{d,n}'$.  However, not all of them contribute to the homology because the existence of a tree $T'<T$ with $T 
\setminus \{x\}\subset T'$ does not guarantee that $T'\in T_{d,n}'$.

Let $T$ be a $3$-valent tree and $V$ be an internal vertex, such that  
the two children of $V$ are the end $m_0$ and a bounded edge leading  
to a vertex $W$ which has the end $m_1$ as child. Assume $m_0>m_1$  
and $m_0$ and $m_1$ are either both marked ends or unmarked ends.
Denote by $T_1$ the subtree of $T$ that can be reached by the parent  
of $V$ and by $T_2$ the subtree that can be reached from $W$ via the  
second child.

\begin{center}
\includegraphics{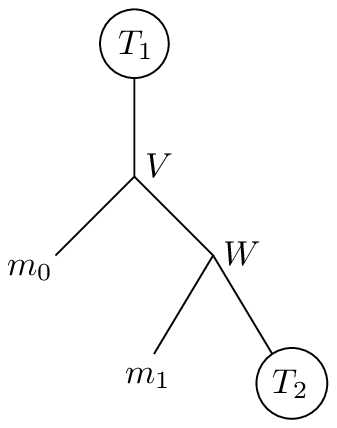}
\end{center}

There are two other trees that differ from $T$ exactly at the split corresponding to edge $\{V,W\}$:
\begin{center}
\includegraphics{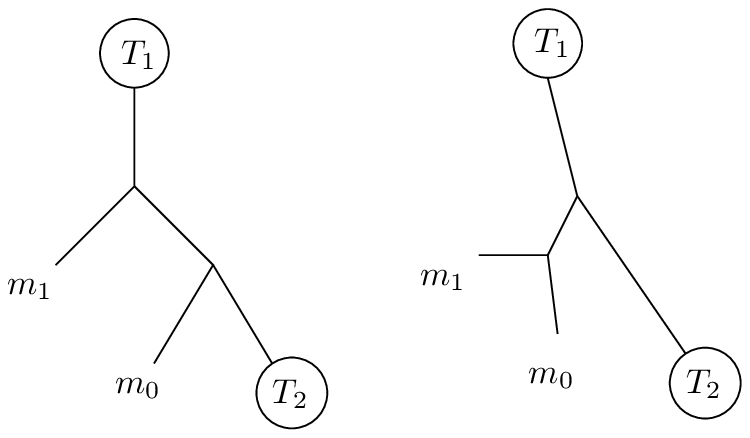}
\end{center}
The first tree is bigger than $T$ and the second is not in $T_{d,n}'$, hence $T$ does not close a loop and does not count toward the top homology of $T_{d,n}'$.

The combs which remain are those which satisfy: if a marked (unmarked) end $m_l$ to the left is followed by a marked (unmarked) end $m_{l+1}$, then $m_l<m_{l+1}$.

We have to show that for those combs we can find $T'$ for each edge  
label $x$.
There are two cases: either $x$ is an edge between two marked (unmarked) ends $m_l$ and $m_{l+1}$ with $m_l<m_{1+1}$, or $m_l$ is  
marked and $m_{1+1}$ unmarked (respectively, the other way round).
In the first case, we just swap $m_l$ and $m_{l+1}$.
\begin{center}
\includegraphics{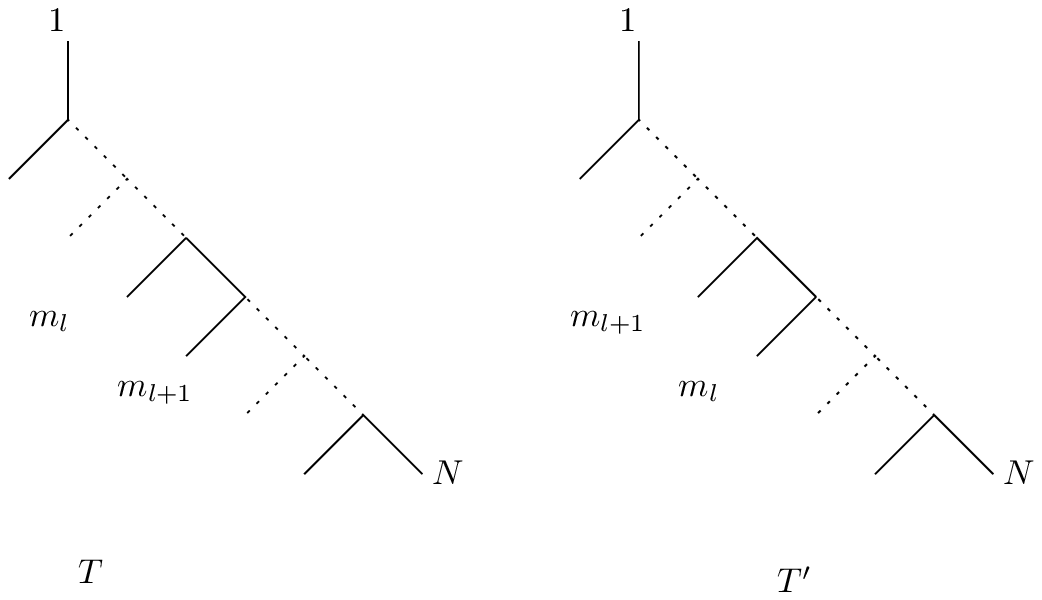}
\end{center}

The new tree $T'$ is still in $T_{d,n}'$. It contains all edge labels  
of $T$ besides the edge label $x=\{m_{l+1},\ldots,m_{n+d-1},N\}$. Instead of $x$, it  
has the edge label $$\{m_l,m_{l+2},\ldots,m_{n+d-1},N\}.$$ Furthermore,  
it is smaller than $T$, because at the vertex above $x$ the right  
subtree contains $m_{l+1}$ which is bigger than $m_l$.

In the second case, we bring $m_l$ and $m_{l+1}$ together to one vertex:
\begin{center}
\includegraphics{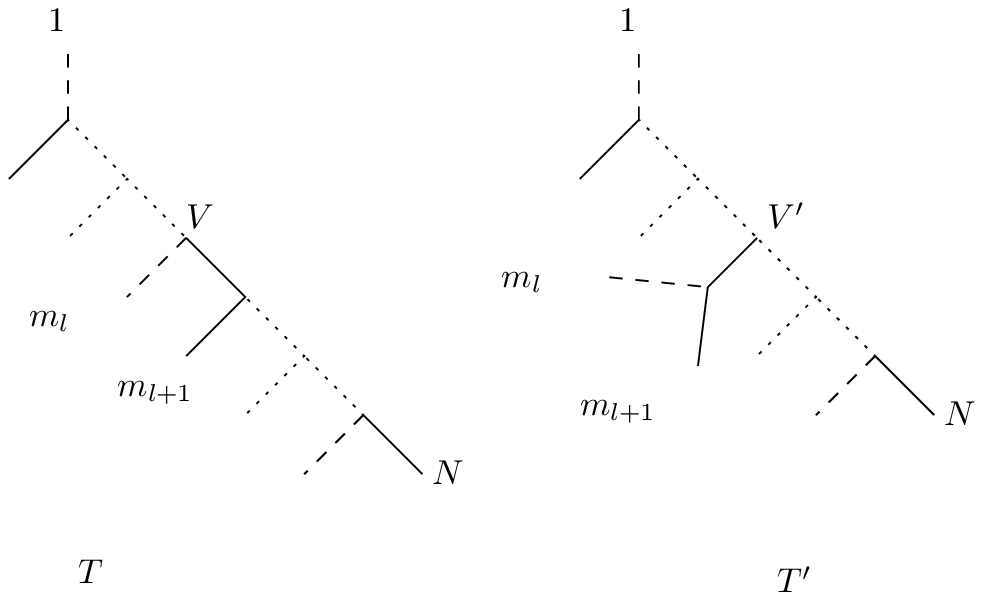}
\end{center}

The edge label $x=\{m_{l+1},\ldots,m_{n+d-2},N\}$ is replaced by $\{m_l,m_{l+1}\}$. All other edge labels remain. Of  
course, $T'\in T_{d,n}'$. Also, $T'<T$, because the right child of the  
vertex $V'$ is smaller than the right child of the vertex $V$, and  
above those vertices, the trees coincide.

Altogether, this proves the following theorem:
\begin{theorem}[Homology of $T_{d,n}''$]\label{thm-hom}
The top homology of $T_{d,n}''$ is $\Z^{h}$, where $h$ is equal to the  
number of combs starting with $1$ and ending with $N$, the neighbour  
of $1$ an unmarked end, the neighbour of $N$ a marked end, two  
consecutive ends either have different colors or have the same color with increasing labels.
\end{theorem}

\begin{corollary}
The rank of the top homology of $T_{d,n}''$ is 
$$
\sum_{k=1}^{\text{min}(n-1,d-1)} \left( \sum_{i=1}^k (-1)^{k-i} {k \choose i} i^{n-1} \right) \left( \sum_{j=1}^k (-1)^{k-j} {k \choose j} j^{d-1} \right).
$$
\end{corollary}

\begin{proof}
It follows from the previous lemma that the rank of the top homology is the number of ways to simultaneously partition an $n-1$-element set and a $d-1$-element set into the same number of nonempty ordered parts.  The number of partitions of an $m$-element set into exactly $k$ nonempty ordered part is 
$$
k! \cdot S(m,k) = \sum_{i=1}^k (-1)^{k-i} {k \choose i} i^{m} 
$$
where $S(m,k)$ are the Stirling numbers of the second kind.
\end{proof}

The computations in \cite{Dev05} for the top homology of $T_{3,n}'', T_{4,4}''$, and $T_{4,5}''$, which are $2^n-3, 73$, and $301$ respectively, agree with our formula.







\end{document}